\documentclass[11pt]{article}
\usepackage{pb-diagram}
\usepackage{amsmath,amsthm,amsfonts,amssymb,latexsym}
\usepackage{amsmath}

\usepackage{a4wide}

\newtheorem{thm}{Theorem}[section]
\newtheorem{cor}[thm]{Corollary}
\newtheorem{lem}[thm]{Lemma}
\newtheorem{prop}[thm]{Proposition}

\newtheorem{rmk}[thm]{Remark}

\numberwithin{equation}{section}

\newtheorem{theorem}{Theorem}[section]

\newtheorem{proposition}[theorem]{Proposition}

\newcommand{\fL}{{\mathfrak{L}}}

\newcommand{\cK}{{\mathcal{K}}}

\newcommand{\Hom}{\mathrm{Hom}}

\newcommand{\Vect}{\mathrm{Vect}}

\newcommand{\half}{\frac{1}{2}}

\def\a{\alpha}

\def\l{\lambda}

\begin{document}


\title{ Differential Operators on  the Weighted Densities on
the Supercircle $S^{1|n}$}


\author{Nader Belghith\thanks{University of Sousse, ISSAT, Sousse, Tunisia. E-mail:
nader.belghith@issatso.rnu.tn} \and Mabrouk Ben
Ammar\thanks{University of Sfax, Faculty  of Sciences, BP 802, 3038
Sfax, Tunisia. E-mail: mabrouk.benammar@fss.rnu.tn}\and Nizar Ben
Fraj\thanks{University of Carthage, IPEIN, Tunisia. E-mail:
benfraj\_nizar@yahoo.fr}}

\maketitle

\begin{abstract}
Over the $(1,n)$-dimensional real supercircle, we consider the
$\mathcal{K}(n)$-modules $\frak{D}^{n,k}_{\lambda,\mu}$ of linear differential operators of order $k$
 acting on the superspaces of weighted
densities, where $\mathcal{K}(n)$ is the Lie superalgebra of contact
vector fields. We give, in contrast to the classical setting, a
classification of these modules for $n=1$. We also prove that
$\frak{D}^{n,k}_{\lambda,\mu}$ and $\frak{D}_{\rho,\nu}^{n,k}$ are
isomorphic for $\rho=\frac{2-n}{2}-\mu$ and
$\nu=\frac{2-n}{2}-\lambda$. This work is the simplest superization
of a result by Gargoubi and Ovsienko [Modules of Differential
Operators on the Real Line, Functional Analysis and Its
Applications, Vol. 35, No. 1, pp. 13--18, 2001.]
\end{abstract}

\bigskip
\thispagestyle{empty}

\section{Introduction.}

Let $\Vect(S^1)$ be the Lie algebra of vector fields on $S^1$.
Consider the 1-parameter deformation of the natural
 $\Vect(S^1)$-action  on $C^\infty(S^1)$:
\begin{equation*}
\mathrm{L}_{X_h}^\lambda(f)=\mathrm{L}_{X_h}(f)+\lambda
h'f,\quad\text{where}\quad X_h =h\frac{d}{dx}.
\end{equation*}
The $\Vect(S^1)$-module so defined is the space ${\cal F}_\lambda$ of weighted densities of
weight $\l\in\mathbb{R}$:
$$
{\cal F}_\lambda=\left\{fdx^{\lambda}\mid
f\in C^\infty(S^1)\right\}.
$$
We denote $\mathrm{D}_{\lambda,\mu}$ the space of linear differential operators from ${\cal F}_\lambda$ to ${\cal F}_\mu$ and $\mathrm{D}_{\lambda,\mu}^k$ the space
of linear differential operators of order $k$ which are naturally endowed with a
$\Vect(S^1)$-module structure defined by the action $\mathrm{L}^{\lambda,\mu}$:
$$
\mathrm{L}_{X_h}^{\lambda,\mu}(A):=\mathrm{L}_{X_h}^\mu\circ
A-A\circ\mathrm{L}_{X_h}^\lambda,\quad\text{ where }\quad A\in\mathrm{D}_{\lambda,\mu}.
$$
Gargoubi and Ovsienko \cite{go} classified these modules and gave a
complete list of isomorphisms between distinct modules
$\mathrm{D}_{\lambda,\mu}^k$. The classification problem for modules
of differential operators on a smooth manifold was posed for
$\lambda=\mu$ and solved for modules of second-order operators in
\cite{dv}. The modules $\mathrm{D}_{\lambda,\lambda}^k$ on
$\mathbb{R}$ were classified in \cite{gov}. In the multidimensional
case, this classification problem was solved in \cite{LMT, M}.

\noindent In this paper we study the simplest super analog of the
problem solved in \cite{go}, namely, we consider the supercircle
$S^{1|n}$ equipped with the contact structure determined by a 1-form
$\alpha_n$, and the Lie superalgebra $\mathcal{K}(n)$ of contact
vector fields on $S^{1|n}$. We introduce the $\mathcal{K}(n)$-module
$\mathfrak{F}^n_\lambda$ of $\lambda$-densities on $S^{1|n}$ and the
$\mathcal{K}(n)$-module of linear differential operators,
$\frak{D}_{\lambda,\mu}^n
:=\Hom_{\rm{diff}}(\mathfrak{F}^n_{\lambda},\mathfrak{F}^n_{\mu})$,
which are super analogs of the spaces $\mathcal{F}_\lambda$ and
$\mathrm{D}_{\lambda,\mu}$ respectively. The $\mathcal{K}(n)$-module
$\frak{D}_{\lambda,\mu}^n$ is filtered:
\begin{equation*}
\mathfrak{D}^{n,0}_{\lambda,\mu}\subset\mathfrak{D}^{n,\frac{1}{2}}_{\lambda,\mu}\subset
\mathfrak{D}^{n,1}_{\lambda,\mu}\subset
\mathfrak{D}^{n,\frac{3}{2}}_{\lambda,\mu}
\subset\cdots\subset\mathfrak{D}^{n,k-\frac{1}{2}}_{\lambda,\mu}\subset\mathfrak{D}^{n,k}_{\lambda,\mu}
\cdots
\end{equation*}
For $n=1$, we omit the subscript $n$, that is, $\mathfrak{F}^n_{\lambda}$ and $\mathfrak{D}^{n,k}_{\lambda,\mu}$ will be simply denoted $\mathfrak{F}_{\lambda}$ and $\mathfrak{D}^{k}_{\lambda,\mu}$.

The aim of this paper is to classify these modules $\mathfrak{D}^{n,k}_{\lambda,\mu}$. For $n=1$ we shall
give a complete list of isomorphisms between distinct modules
$\mathfrak{D}^{k}_{\lambda,\mu}$. Moreover, we prove that
$\mathfrak{D}_{\lambda,\mu}^{n,k}$ and
$\mathfrak{D}_{\rho,\nu}^{n,k}$ are $\cK(n)$-isomorphic for
$\nu=\frac{2-n}{2}-\lambda$ and $\rho=\frac{2-n}{2}-\mu$. The
complete classification of modules
$\mathfrak{D}^{n,k}_{\lambda,\mu}$, for $n\geq2$, needs an other
study. We mention that a similar problem was considered in \cite{c1}
for the case of pseudodifferential operators instead of differential
operators.

\section{The main definitions.}
In  this section, we recall the main definitions and facts related
to the geometry of the supercircle $S^{1\mid 1}$; for more details,
see \cite{c, gmo, gls, L}.
\subsection{The Lie superalgebra of contact vector fields on
$S^{1|1}$}
Let $S^{1\mid 1}$ be the supercircle with local coordinates
$(x,\theta),$ where  $\theta$ is an odd indeterminate: $\theta^2=0$.
We introduce the vector fields
${\eta}=\partial_{\theta}+\theta\partial_x$ and
$\overline{\eta}=\partial_{\theta}-\theta\partial_x$. The
supercircle $S^{1\mid 1}$ is equipped with the standard contact
structure given by the distribution $\langle\overline{\eta}\rangle$.
That is, the distribution $\langle\overline{\eta}\rangle$ is the
kernel of the following $1$-form:
\begin{equation*}
\a=dx+\theta d\theta.
\end{equation*}
On $C^\infty(S^{1|1})$, we consider the contact bracket
\begin{equation*}
\{F,G\}=FG'-F'G-\frac{1}{2}(-1)^{|F|}\overline{\eta}(F)
\overline{\eta}(G),
\end{equation*}
where the subscript $'$ stands for $\frac{d}{dx}$ and $|F|$ is the parity of an homogeneous function
$F$.
Let $\mathrm{Vect}(S^{1|1})$ be the superspace of
 vector fields on $S^{1|1}$:
\begin{equation*}\mathrm{Vect}(S^{1|1})=\left\{F_0\partial_x
+  F_1\partial_{\theta} \mid ~F_0,\,F_1\in C^\infty(S^{1|1})
\right\},\end{equation*}  and consider the superspace
$\mathcal{K}(1)$ of contact vector fields on $S^{1|1}$ (also known
as the Neveu-Schwartz superalgebra without central charge, cf.
\cite{AB, R}). That is, $\mathcal{K}(1)$ is the superspace of vector
fields on $S^{1|1}$ preserving the distribution
$\langle\overline{\eta}\rangle$:
$$
\mathcal{K}(1)=\big\{X\in\mathrm{Vect}(S^{1|1})~|~[X,\,\overline{\eta}]=
F_X\overline{\eta}\quad\hbox{for some}~F_X\in
C^\infty(S^{1|1})\big\}. $$ The Lie superalgebra $\mathcal{K}(1)$ is
spanned by the vector fields of the form:
\begin{equation*}
\label{field} X_F=F\partial_x
-\frac{1}{2}(-1)^{|F|}\overline{\eta}(F)\overline{\eta},\quad\text{where}\quad
F\in C^\infty(S^{1|1}).
\end{equation*}
Of course, $\mathcal{K}(1)$ is a subalgebra of
$\mathrm{Vect}(S^{1|1})$, and $\mathcal{K}(1)$ acts on
$C^{\infty}(S^{1|1})$ through:
\begin{equation*}
\mathfrak{L}_{X_F}(G)=FG'-\frac{1}{2}(-1)^{|F|}\overline{\eta}(F)\cdot\overline{\eta}(G).
\end{equation*}
The bracket in $\mathcal{K}(1)$ can be written as: $
[X_F,\,X_G]=X_{\{F,\,G\}}$.
\subsection{The space of weighted densities on
$S^{1|1}$}

In the super setting, by replacing $dx$ by the 1-form $\alpha$, we get analogous definition for
weighted densities {\sl i.e.} we define the space of $\lambda$-densities as
\begin{equation*}
\mathfrak{F}_\lambda=\left\{F(x,\theta)\alpha^\lambda~~|~~F(x,\theta) \in C^\infty(S^{1|1})\right\}.
\end{equation*}
As a vector space, $\mathfrak{F}_\lambda$ is isomorphic to $C^\infty(S^{1|1})$,
but the Lie derivative of the density $G\alpha^\lambda$ along the vector field $X_F$ in
$\mathcal{K}(1)$ is:
\begin{equation*}
\mathfrak{L}^{\lambda}_{X_F}(G\alpha^\lambda)=\left(\mathfrak{L}_{X_F}(G)+ \lambda F'G\right)\alpha^\lambda.
\end{equation*}
 Obviously $\mathcal{K}(1)$ and $\mathfrak{F}_{-1}$ are isomorphic as $\mathcal{K}(1)$-modules.
 Naturally $\cK(1)$ acts on the superspace
$\frak{D}_{\lambda,\mu}
:=\Hom_{\rm{diff}}(\mathfrak{F}_{\lambda},\mathfrak{F}_{\mu})$ through:
\begin{equation}\label{d-action}
\fL^{\l,\mu}_{X_F}(A)=\fL^{\mu}_{X_F}\circ A-(-1)^{|A||F|} A\circ
\fL^{\l}_{X_F}.
\end{equation}
Since $\overline{\eta}^2=-\partial_x$, any differential operator
$A\in\mathfrak{D}_{\lambda,\mu}$ can be expressed in the form:
\begin{equation}\label{aaa}
A(F\alpha^\lambda)=\sum_{i=0}^\ell a_i\,
\overline{\eta}^i(F)\alpha^\mu,
\end{equation}
where the coefficients $a_i\in C^\infty(S^{1|1})$ and
$\ell\in\mathbb{N}$.
For $k\in \frac{1}{2}\mathbb{N}$, the space of differential
operators of the form \eqref{aaa} with $\ell=2k$ is denoted by
$\mathfrak{D}^{k}_{\lambda,\mu}$ and called the space of
differential operators of order $k$.
Thus, we have a $\mathcal{K}(1)$-invariant filtration:
\begin{equation}\label{graded}
\mathfrak{D}^{0}_{\lambda,\mu}\subset\mathfrak{D}^{\frac{1}{2}}_{\lambda,\mu}\subset
\mathfrak{D}^{1}_{\lambda,\mu}\subset\mathfrak{D}^{\frac{3}{2}}_{\lambda,\mu}
\subset\cdots\subset\mathfrak{D}^{i-\frac{1}{2}}_{\lambda,\mu}\subset\mathfrak{D}^{i}_{\lambda,\mu}
\cdots
\end{equation}
The quotient module
$\mathfrak{D}^{i}_{\lambda,\mu}/\mathfrak{D}^{i-\frac{1}{2}}_{\lambda,\mu}$ is
isomorphic to the module of weighted densities
$\Pi^{2i}\left(\mathfrak{F}_{\mu-\lambda-i}\right)$ (see, e.g., \cite{gmo}), where $\Pi$ is the change of parity map. Thus, the graded $\mathcal{K}(1)$-module
$\mathrm{gr}\mathfrak{D}_{\lambda,\mu}$ associated with the
filtration (\ref{graded}) is a direct sum of density modules:
\[
\mathrm{gr}\mathfrak{D}_{\lambda,\mu}=\bigoplus_{i=0}^\infty\Pi^i\left(\mathfrak{F}_{\mu-\lambda-\frac{i}{2}}\right).
\]
Note that this module depends only on the shift, $\mu-\lambda$, of the weights and not on $\mu$ and
$\lambda$ independently. We call this $\mathcal{K}(1)$-module the space of symbols of differential operators and
denote it $\mathfrak{S}_{\mu-\lambda}$. 
The space of symbols of order $\leq k$ is
\begin{equation*}
{\frak
S}^k_{\mu-\lambda}=\bigoplus_{i=0}^{2k}\Pi^i\left(\mathfrak{F}_{\mu-\lambda-\frac{i}{2}}\right).
\end{equation*}
\section{Classification results}
We now give a complete classification of the modules
$\mathfrak{D}^{k}_{\lambda,\mu}$. First note that the difference
$\delta = \mu - \lambda$ of weight is an invariant: the condition
$\mathfrak{D}^{k}_{\lambda,\mu}\simeq\mathfrak{D}^{k}_{\rho,\nu}$
implies that $\mu - \lambda=\nu - \rho$. This is a consequence of
the equivariance with respect to the vector field $X_x$. Moreover,
recall that, for every $k\in\half \mathbb{N}$, there exists a
$\mathcal{K}(1)$-invariant conjugate map from
$\mathfrak{D}_{\lambda,\mu}^{k}$ to
$\mathfrak{D}_{\frac{1}{2}-\mu,\frac{1}{2}-\lambda}^{k}$ defined by:
\[
a\overline{\eta}^i\mapsto(-1)^{[\frac{i+1}{2}]+i|a|}\overline{\eta}^i\circ
a.
\]
Clearly, this map is a $\mathcal{K}(1)$-isomorphism. The module
$\mathfrak{D}_{\frac{1}{2}-\mu,\frac{1}{2}-\lambda}^{k}$ is called
the adjoint module of $\mathfrak{D}^k_{\l,\mu}$. A module with
$\lambda+\mu=\half$ is said to be self-adjoint. We say that a
modules $\mathfrak{D}_{\lambda,\mu}^{k}$ is singular if it is only
isomorphic to its adjoint module.

\medskip

Our main result of this paper is the following:
\begin{theorem}
\label{main11}
\begin{itemize}
\item [i)] For $k\leq2$, the $\mathcal{K}(1)$-modules
$\mathfrak{D}_{\lambda,\mu}^{k}$ and $\mathfrak{D}_{\rho,\nu}^{k}$ are isomorphic if and only if
$ \mu -\lambda=\nu-\rho$ except for the modules listed in the
following table and their adjoint modules
 which are all singular.

\vskip.3cm\centerline {Table 1} \vskip.3cm
\centerline{\begin{tabular}{|c|c|c|c|c|}
            \hline
            $k$ &  $\frac{1}{2}$&1&$\frac{3}{2}$&2 \\[3pt]\hline
            $(\lambda,\,\mu)$ & $(0,\frac{1}{2})$&$(0,\frac{1}{2})$&$ (0,\mu), (-\frac{1}{2}, 1)$
            &$ (0,\mu), (\lambda,\frac{1}{2}-\lambda), (\lambda, \lambda+2)$\\[3pt]\hline
\end{tabular}}
\item [ii)]For $k > 2$, the $\mathcal{K}(1)$-modules
$\mathfrak{D}_{\lambda,\mu}^{k}$ are all singular. \end{itemize}
\end{theorem}
The proof of Theorem \ref{main11} will be the subject of  sections
\ref{proof1} and \ref{proof2}. In fact, we need first to study the
action of $\mathcal{K}(1)$ on
$\mathfrak{D}_{\lambda,\mu}^k$
in terms of $\mathfrak{osp}(1|2)$-equivariant symbols.
\section{Modules of differential operators over $\mathfrak{osp}(1|2)$}
Consider the Lie superalgebra $\mathfrak{osp}(1|2)\subset
\mathcal{K}(1)$ generated by the functions: $1,\, x,\, x^2,\,
\theta,\, x \theta$. The Lie superalgebra $\mathfrak{osp}(1|2)$
plays a special role and allows one to identify
$\mathfrak{D}_{\lambda,\mu}^k$ with ${\mathfrak S}_{ \mu - \lambda}^k$, in a canonical way. The following
result (see \cite{gmo}) shows that, for generic values of $\lambda$ and $\mu$,
$\mathfrak{D}_{\lambda,\mu}^k$ and ${\mathfrak{S}}_{ \mu - \lambda}^k$ are
isomorphic as $\mathfrak{osp}(1|2)$-modules.
\begin{thm}\cite{gmo}.
\label{gmo1} (i) If $ \mu - \lambda$ is nonresonant, i.e., $ \mu - \lambda
\notin\frac{1}{2}\mathbb{N }\setminus\{0\}$, then
$\mathfrak{D}_{\lambda,\mu}$ and $\mathfrak{S}_{ \mu - \lambda}$ are $\mathfrak{osp}(1|2)$-isomorphic through the
unique $\mathfrak{osp}(1|2)$-invariant symbol map $\sigma_{\lambda,\mu}$ defined by:
\begin{equation}
\label{valentin} \displaystyle
\sigma_{\lambda,\mu}(a\overline{\eta}^k) = \sum_{n=0}^k \Pi^{k-n}\left(\gamma_n^k
\eta^n(a) \alpha^{\mu-\lambda - \frac{k-n}{2}}\right),
\end{equation}
where
\begin{equation}
\label{valentain}
\gamma_n^k=(-1)^{\lbrack \frac{n+1}{2}\rbrack}\frac{\left(\begin{array}{c}\lbrack \frac{k}{2}\rbrack\\
\lbrack \frac{2n+1-(-1)^{n+k}}{4}\rbrack \end{array}\right)
\left(\begin{array}{c}\lbrack \frac{k-1}{2}\rbrack+2\lambda\\
\lbrack \frac{2n+1+(-1)^{n+k}}{4}\rbrack \end{array}\right)}{\left(\begin{array}{c}2(\mu - \lambda)+n-k-1\\
\lbrack \frac{n+1}{2}\rbrack \end{array}\right)}
\end{equation}
with $\big(^\nu_i\big)=\frac {\nu(\nu-1)\cdots (\nu-i+1)}{i!}$ and
$[x]$ denotes the integer part of a real number $x$.

\medskip

(ii) In the resonant cases the $\mathfrak{osp}(1|2)$-modules
$\mathfrak{D}_{\lambda,\mu}$ and $\mathfrak{S}_{ \mu - \lambda}$ are not
isomorphic, except for $(\lambda, \mu) = (\frac{1 - m}{4},
\frac{1+m}{4})$, where $m$ is an odd integer.
\end{thm}
The main idea of proof of Theorem~\ref{main11} is to use the
$\mathfrak{osp}(1|2)$-equivariant symbol mapping $\sigma_{\l,\mu}$
to reduce the action of $\mathcal{K}(1)$ on
$\mathfrak{D}_{\lambda,\mu}^k$ to a canonical form. In other words,
we shall use the following diagram
\[
\begin{diagram}
\node{\mathfrak{D}^k_{\lambda,\mu}}
\arrow[4]{e,t}{{\fL^{\l,\mu}_{X_F}}}
\arrow[2]{s,l}{{\sigma}_{\lambda,\mu}}
\node[4]{\mathfrak{D}^k_{\lambda,\mu}}
   \arrow[2]{s,r}{{\sigma}_{\lambda,\mu}}
\\
\\
\node{\mathfrak{S}^k_{ \mu - \lambda}}
\arrow[4]{e,t}{{\widetilde{\fL}^{\l,\mu}_{X_F}}}
   \node[4]{\mathfrak{S}^k_{ \mu - \lambda}}
\end{diagram}
\]
and compare the action
$\widetilde{\fL}^{\l,\mu}_{X_F}:=\sigma_{\lambda,\mu} \circ
\fL^{\l,\mu}_{X_F} \circ \sigma_{\lambda,\mu}^{-1}$ with the
standard action of $\mathcal{K}(1)$ on $\mathfrak{S}_{ \mu - \lambda}^k$.
\section{The action of $\mathcal{K}(1)$ in the $\mathfrak{osp}(1|2)$-invariant form.}
The action of $\mathcal{K}(1)$ on
$\mathfrak{D}_{\lambda,\mu}^k$ 
in terms of $\mathfrak{osp}(1|2)$-equivariant symbols is closely
related to the space $\mathcal{S}_k^\lambda$ of
$\mathfrak{osp}(1|2)$-invariant linear operators from $\cK{(1)}$ to
$\frak{D}_{\lambda,\lambda+k-1}$ vanishing on $\mathfrak{osp}(1|2)$.
For $k>2$, the space $\mathcal{S}_k^\lambda$ is one dimensional,
spanned by the maps:
$$X_F\mapsto\left(G\alpha^\lambda\mapsto\frak{J}_{k}^{\lambda}(F,G)\alpha^{\lambda+k-1}\right)$$
where $\frak{J}_{k}^{\lambda}$ is
the supertransvectant $\frak{J}_{k}^{-1,\lambda}$ defined in \cite{bbbbk} (see also \cite{go1, gt}).
The operators $\frak{J}_{k}^{\lambda}$ labeled by semi-integer $k$
are odd and they are given by:
\begin{equation*}\begin{array}{lllll}
\frak{J}_{k}^{\lambda}(F,G)&=\displaystyle\sum_{i+j=[k],i\geq2}
\Gamma_{i,j,k}^{\lambda}\left((-1)^{|F|}([k]-j-2)F^{(i)}\overline{\eta}(G^{(j)})-
(2\l+[k]-i)\overline{\eta}(F^{(i)})G^{(j)}\right).
\end{array}\end{equation*}
The operators $\frak{J}_{k}^{\lambda}$, where $k\in \mathbb{N}$, are
even and they are given by:
\begin{equation*}\begin{array}{lllll}
\frak{J}_k^{\lambda}(F,G)&=\displaystyle\sum_{i+j=k-1,
i\geq2}(-1)^{|F|}
\Gamma_{i,j,k-1}^{\lambda}\overline{\eta}(F^{(i)})\overline{\eta}(G^{(j)})
-\displaystyle\sum_{i+j=k,i\geq3}\Gamma_{i,j,k-1}^{\lambda}F^{(i)}G^{(j)},\end{array}
\end{equation*} where $\big(^x_i\big)=\frac {x(x-1)\cdots (x-i+1)}{i!}$ and $[k]$ denotes
the integer part of $k$, $k>0$, and
\begin{equation*}
\label{coe} \Gamma_{i,j,k}^{\lambda}=(-1)^{j}
\begin{pmatrix}[k]-2\\j\end{pmatrix}
\begin{pmatrix}2\l+[k]\\i\end{pmatrix}.
\end{equation*}
We will need the expressions of $\frak{J}_{\frac{5}{2}}^{\lambda}$,
$\frak{J}_{3}^{\lambda}$ and $\frak{J}_{\frac{7}{2}}^{\lambda}$:
\begin{equation}
\small{ \left\{\begin{array}{lllllllll}
\frak{J}_{\frac{5}{2}}^{\lambda}(F,G)&=&\overline{\eta}(F'')G,
\\[2pt]
\frak{J}_{3}^{\lambda}(F,G)&=&\left({2\over3}\lambda
F^{(3)}G-(-1)^{|F|}\overline{\eta}(F'') \overline{\eta}(G)
\right),
\\[2pt]
\frak{J}_{\frac{7}{2}}^{\lambda}(F,G)&=&\left(
2\lambda\overline{\eta}(F^{(3)})G-3\overline{\eta}(F'')G'
-(-1)^{|F|}F^{(3)}\overline{\eta}(G) \right),
\end{array}\right.}
\end{equation}

Now, we compute  the action of $\mathcal{K}(1)$ on
$\mathfrak{D}_{\lambda,\mu}^k$ in terms of
$\mathfrak{osp}(1|2)$-equivariant symbols. First, using formula (\ref{d-action}) and the graded
Leibniz formula:
\begin{equation}
\label{libniz} \overline{\eta}^j\circ
F=\sum_{i=0}^j\begin{pmatrix}j\\i\end{pmatrix}_s(-1)^{|F|(j-i)}\overline{\eta}^i(F)\overline{\eta}^{j-i},
\quad\text{where}\quad\begin{pmatrix}j\\i\end{pmatrix}_s=
\left\{
\begin{array}{llllllll}
{\begin{pmatrix}\,[{j\over2}]\\\,[{i\over2}]\end{pmatrix}}
&\text{if}~ i~ \hbox{is even
or}~ j~ \hbox{is odd},\\[10pt]
~~~\,0 &\hbox{otherwise}.
\end{array}
\right.
\end{equation}
we easily check the following result:
\begin{lem}
\label{a}
The natural action of  $\mathcal{K}(1)$ on
${\frak{D}}_{\lambda,\mu}^k$ is given by $
\fL^{\l,\mu}_{X_F}(A):=\sum_{i=0}^{2k} a^{X_F}_i\,
\overline{\eta}^i,$ where
\begin{equation}
\label{naturel action}
\begin{array}{lllllll}
a^{X_F}_i&=&\fL^{\mu-\lambda-{i\over 2}}_{X_F}(a_i)-
\displaystyle\sum_{j\geq
i+1}^{2k}(-1)^{(|F|+|a_j|)(j-i)}\zeta_{i,j,\lambda}\overline{\eta}^{j-i}(F')a_j\end{array}
\end{equation} with
\begin{equation}
\label{zeta} {\small\small
\zeta_{i,j,\lambda}=\lambda{\begin{pmatrix}j\\j-i\end{pmatrix}}_s
-{(-1)^{i}\over 2}\begin{pmatrix}j\\j-i+1\end{pmatrix}_s
+{\begin{pmatrix}j\\j-i+2\end{pmatrix}}_s}.
\end{equation}
\end{lem}

Now, we need to study the action of $\mathcal{K}(1)$ over
$\mathfrak{D}_{\lambda,\mu}^k$ in terms of
$\mathfrak{osp}(1|2)$-equivariant symbols, thus, let
\begin{equation}
\label{ipem}
\widetilde{\fL}^{\l,\mu}_{X_F}\left(\sum_{p=0}^{2k}\Pi^{p}\left(P_p\alpha^{\delta-{p\over2}}\right)\right)
=\sum_{p=0}^{2k}\Pi^{p}\left(P^{X_F}_p\alpha^{\delta-{p\over2}}\right),\quad\text{where}\quad \delta=\mu-\lambda,
\end{equation}
then, we need to compute the terms $P^{X_F}_p$.
\begin{prop}
\label{nizar0} (i) The terms $P^{X_F}_p$ are given by:
\begin{equation}
\label{bellie} P_p^{X_F} = \fL^{\delta - \frac{p}{2}}_{X_F} (P_p) +
\sum_{j=p+3}^{2k}\beta_p^j\,
\pi^{j-p}\circ\frak{J}_{\frac{j-p}{2}+1}^{\delta - \frac{j}{2}
}\left(X_F,P_j\right),
\end{equation}
where
$\pi(F)=(-1)^{|F|}F$ and the coefficients $\beta_p^j$ are some functions of $\lambda$ and $\mu$.

\medskip

(ii) For $j\leq5$, the coefficients $\beta_p^j$ of formula
(\ref{bellie}) are given by:
\begin{equation}
\begin{array}{lll}
\label{coffbeta5} \beta_0^3
= \beta_0^3({\lambda, \mu}) =-
\frac{\lambda(2\delta+2\lambda-1)}{2\delta-2},
\\[3pt]
\beta_0^4 = \beta_0^4({\lambda, \mu}) =
-\frac{3\lambda(2\delta+2\lambda-1)}{(2\delta-1)(2\delta-4)},
\\[3pt]
\beta_1^4 =  \beta_1^4({\lambda, \mu}) =-
\frac{2\delta+4\lambda-1}{2(2\delta-3)},
\\[3pt]
\beta_0^5 = \beta_0^5({\lambda, \mu}) =
\frac{\lambda(2\delta+2\lambda-1)(2\delta+4\lambda-1)}{(2\delta-1)(2\delta-3)(2\delta-5)},
\\[3pt]
\beta_1^5 = \beta_1^5({\lambda, \mu}) = -\frac{3(4\lambda \delta + 2
\delta +4\lambda^2-2\lambda-1)}{(2\delta-5)(4\delta-4)},
\\[3pt]
\beta_2^5 =  \beta_2^5({\lambda, \mu}) =- \frac{\delta+4\lambda
\delta-2\lambda+4\lambda^2}{2(\delta-2)}.
\end{array}
\end{equation}

\end{prop}

\begin{proofname}.
(i) According to Lemma \ref{a} and formula
\eqref{valentin}, we prove that $P_p^{X_F}$ can be expressed as follows
\begin{equation}\label{trivial}
P_p^{X_F} = \fL^{\delta - \frac{p}{2}}_{X_F} (P_p) + \sum_{j=0}^{2k}f_j\left(X_F,P_j\right),
\end{equation}
where the $f_j$ are bilinear maps from $\mathcal{K}(1)\times
\mathfrak{F}_{\delta-\frac{j}{2}}$ to
$\mathfrak{F}_{\delta-\frac{p}{2}}$ vanishing on
$\mathfrak{osp}(1|2)$. Since $\widetilde{\fL}^{\l,\mu}$ is a
$\mathcal{K}(1)$-action on $\mathfrak{S}_\delta^k$, then $f_j$ has
the same parity as the integer $j-p$. Moreover, for $X_G\in
\mathfrak{osp}(1|2)$, we have
\[
\mathfrak{L}_{X_G}\circ\widetilde{\mathfrak{L}}_{X_F}^{\lambda,\mu}=
\widetilde{\mathfrak{L}}_{[X_G,X_F]}^{\lambda,\mu}+ (-1)^{|F||G|}\widetilde{\mathfrak{L}}_{X_F}^{\lambda,\mu}\circ\mathfrak{L}_{X_G}.
\]
Thus, from \eqref{trivial}, we deduce that
\[
\mathfrak{L}_{X_G}^{\delta-\frac{p}{2}}f_j(X_F,P_j)=f_j([X_G,X_F],P_j)+(-1)^{|F||G|} f_j(X_F,\mathfrak{L}_{X_G}^{\delta-\frac{j}{2}}(P_j)).
\]
Therefore, the map $\pi^{j-p}\circ f_j$ is a supertranvectant vanishing on $\mathfrak{osp}(1|2)$.
Thus, up to a scalar factor, we have  $f_j=\pi^{j-p}\circ\frak{J}_{\frac{j-p}{2}+1}^{\delta - \frac{j}{2}}$
for $j\geq p+3$; otherwise $f_j=0$. 

(ii) By a direct computation, using formulas (\ref{naturel action}) and \eqref{valentin}, we get expression (\ref{coffbeta5}).\hfill$\Box$
\end{proofname}

\medskip
\section{Proof of Theorem \ref{main11} in the generic case}
\label{proof1} In this section we prove  Theorem \ref{main11} for the
nonresonant values of  $\delta=\mu-\lambda$.

\begin{prop}
\label{bellie village} For $k\,\in\,\frac{1}{2}\mathbb{N}$, let $T:
\mathfrak{D}_{\lambda,\mu}^k \rightarrow \mathfrak{D}_{\rho,\nu}^k$
 be an isomorphism of $\mathcal{K}(1)$-modules. Then  the linear
mapping $\sigma_{\rho, \nu} \circ T \circ \sigma_{\lambda,
\mu}^{-1}$ on $\mathfrak{S}_\delta^k$ is diagonal and the
$\Pi^i\left(\mathfrak{F}_{\delta-\frac{i}{2}}\right)$ are
eigenspaces:
{\begin{equation}\label{diag1}\Pi^{i}\left(P_{i}^T\alpha^{\delta-\frac{i}{2}}\right):=
\sigma_{\rho, \nu} \circ T \circ \sigma_{\lambda, \mu}^{-1}\left(
\Pi^{i}\left(P_{i}\alpha^{\delta-\frac{i}{2}}\right)\right) =
\Pi^{i}\left(\tau_{i} P_{i}\alpha^{\delta-\frac{i}{2}}\right),\;
\tau_i\in\,\mathbb{R}\setminus\{0\}.
\end{equation}}
\end{prop}
\begin{proofname}.
Since $T$ is an isomorphism of $\mathcal{K}(1)$-modules, it is also
an isomorphism of $\mathfrak{osp}(1|2)$-modules. The uniqueness of
the $\mathfrak{osp}(1|2)$-equivariant symbols mapping shows that the
linear mapping $\sigma_{\rho, \nu} \circ T \circ \sigma_{\lambda,
\mu}^{-1}$ on $\mathfrak{S}_\delta^k$ is diagonal and the
$\Pi^i\left(\mathfrak{F}_{\delta-\frac{i}{2}}\right)$ are
eigenspaces.\hfill$\Box$
\end{proofname}
\subsection{The construction of isomorphisms}
To prove Theorem~\ref{main11}, we construct the desired isomorphism
explicitly in terms of projectively equivariant symbols using Proposition \ref{nizar0} and Proposition \ref{bellie village}.
\begin{itemize}
  \item [i)]
 For $k = \frac{1}{2}$, formula (\ref{diag1}) defines an
isomorphism $T: \mathfrak{D}_{\lambda,\mu}^{\frac{1}{2}} \rightarrow
\mathfrak{D}_{\rho,\nu}^{\frac{1}{2}}$
for all $\tau_0, \tau_1 \neq 0$:
{\small\begin{equation}
\label{isok1} \left(P_0^T\alpha^\delta,\, \Pi\left(P_1^T\alpha^{\delta-\frac{1}{2}}\right)\right) = \left(\tau_0 P_0\alpha^\delta,\,  \Pi\left(\tau_1 P_1\alpha^{\delta-\frac{1}{2}}\right)\right),
\end{equation}}
since the action  (\ref{bellie}) is
{\small$$
\left({P}_0^{X_F}, {P}_1^{X_F}\right) =
\left(\fL^{\delta}_{X_F}({P}_0),~ \fL^{\delta-\half}_{X_F}({P}_1)\right).
$$}
\item [ii)] For $k = 1$, formula (\ref{diag1}) defines an
isomorphism $T: \mathfrak{D}_{\lambda,\mu}^1 \rightarrow
\mathfrak{D}_{\rho,\nu}^1$
for all $\tau_0, \tau_1, \tau_2 \neq 0$:
{\small\begin{equation}
\label{isok2} \left(P_0^T\alpha^\delta,~ \Pi\left(P_1^T\alpha^{\delta-\frac{1}{2}}\right),~ P_2^T\alpha^{\delta-1}\right) = \left(\tau_0 P_0\alpha^\delta,~ \Pi\left(\tau_1 P_1\alpha^{\delta-\frac{1}{2}}\right),~
\tau_2 P_2\alpha^{\delta-1}\right),
\end{equation}}
since the action  (\ref{bellie}) is
{\small$$\left({P}_0^{X_F}, {P}_1^{X_F},
P_2^{X_F}\right) = \left(\fL^{\delta}_{X_F}({P}_0),~
\fL^{\delta-\half}_{X_F}({P}_1),~ \fL^{\delta - 1}_{X_F}({P}_2)\right).
$$}
\item [iii)] For $k = \frac{3}{2}$, let $T: \mathfrak{D}_{\lambda,\mu}^{\frac{3}{2}}
\rightarrow \mathfrak{D}_{\rho,\nu}^{\frac{3}{2}}$ be any
$\cK(1)$-isomorphism defined by: {\small\begin{equation}
\label{isok3}  \Pi^i\left(P_i^T\alpha^{\delta-\frac{i}{2}}\right) =
\Pi^i\left(\tau_i P_i\alpha^{\delta-\frac{i}{2}}\right),\quad i=0,
1,\,2,\,3,\, \tau_i\in \mathbb{R}\setminus\{0\}.
\end{equation}}
Since the action  (\ref{bellie}) is defined, in this case, by:
{\small\[
\begin{array}{lllll}
&{P}_0^{X_F}=\fL^{\delta}_{X_F}({P}_0) + \beta_0^3 \pi\circ
\frak{J}_{\frac{5}{2}}^{\delta-\frac{3}{2}}(F,{P}_3),\\&{P}_1^{X_F}=\fL^{\delta-\half}_{X_F}({P}_1),\quad
{P}_2^{X_F}=\fL^{\delta-1}_{X_F}({P}_2),\quad
{P}_3^{X_F}=\fL^{\delta-\frac{3}{2}}_{X_F}({P}_3),
\end{array}
\]}
the equivariant conditions of $T$ lead to the following condition:
\[
\tau_0\beta_0^3(\lambda, \mu) = \tau_3\beta_0^3(\rho,\nu).
\]
Thus, we distinguish two cases:\\
\begin{itemize}
    \item [(1)] If $\beta_0^3(\lambda,\mu) \neq 0$, then we get
a family of isomorphisms $T:
\mathfrak{D}_{\lambda,\mu}^{\frac{3}{2}} \rightarrow
\mathfrak{D}_{\rho,\nu}^{\frac{3}{2}}$  given by (\ref{diag1}) with
{\small\begin{equation} \label{isok3} \tau_0,\,\tau_1,\,\tau_2 ,\,
\tau_3\in\mathbb{R}\setminus\{0\}\quad\text{and}\quad \tau_0 =
\tau_3 \frac{\beta_0^3(\rho, \nu)}{\beta_0^3(\lambda, \mu)}.
\end{equation}}
    \item [(2)] If $\beta_0^3(\lambda,\mu) = 0$, that is,  $\lambda=0$ or $\mu=\frac{1}{2}$ then, we have $\beta_0^3(\rho,\nu)=0$, so, the modules $ \mathfrak{D}_{\lambda,\mu}^{\frac{3}{2}}$ and $\mathfrak{D}_{\rho,\nu}^{\frac{3}{2}}$ are equal or conjugate. Thus $\mathfrak{D}_{0,\mu}^{\frac{3}{2}}$ is singular.
  \end{itemize}

\item [iv)] For $k = 2$, we have
{\small\[
\begin{array}{lllll}
{P}_0^{X_F} &=& \fL^{\delta}_{X_F}({P}_0) + \beta_0^3 \pi\circ
\frak{J}_{\frac{5}{2}}^{\delta-\frac{3}{2}}(F,{P}_3)  + \beta_0^4
\frak{J}_{3}^{\delta-2}(F,{P}_4),\quad {P}_2^X =
\fL^{\delta-1}_{X_F}({P}_2)
\\
{P}_1^{X_F} &=& \fL^{\delta-\half}_{X_F}({P}_1) + \beta_1^4 \pi\circ
\frak{J}_{\frac{5}{2}}^{\delta-2}(F,{P}_4),\quad
{P}_3^{X_F} = \fL^{\delta-\frac{3}{2}}_{X_F}({P}_3),\quad
{P}_4^{X_F} = \fL^{\delta-2}_{X_F}({P}_4).
\end{array}
\]}
Thus, we get the following conditions:
\begin{equation}
\label{isok} \tau_0 \beta_0^3(\lambda, \mu)= \tau_3\beta_0^3(\rho,
\nu),\quad \tau_0 \beta_0^4(\lambda, \mu)= \tau_4\beta_0^4(\rho,
\nu),\quad \tau_1\beta_1^4(\lambda, \mu) = \tau_4\beta_1^4(\rho,
\nu)
\end{equation}
and then, as in the previous case, we have to distinguish two cases:
\begin{itemize}
  \item [(1)] If $\beta_0^3,~ \beta_0^4,~ \beta_1^4 \neq
0$, then we get a family of isomorphisms $T:
\mathfrak{D}_{\lambda,\mu}^{2} \rightarrow
\mathfrak{D}_{\rho,\nu}^{2}$  given by (\ref{diag1}) with
{\small\begin{equation*} \label{isok3} \tau_2,\,
\tau_3\in\mathbb{R}\setminus\{0\},\,\tau_4=\tau_3,\quad  \tau_0 =
\tau_3\frac{\beta_0^4(\rho, \nu)}{\beta_0^4(\lambda,
\mu)},\quad\text{and}\quad\tau_1 = \tau_3\frac{\beta_1^4(\rho,
\nu)}{\beta_1^4(\lambda, \mu)}.
\end{equation*}}
  \item [(2)] If $\beta_0^3=0$ or $\beta_0^4=0$ or $\beta_1^4 =0$, then, as in the previous case, we prove
  that the modules $\mathfrak{D}_{0,\mu}^{2}$ and $\mathfrak{D}_{\lambda,\frac{1}{2}-\lambda}^{2}$ are singular.
\end{itemize}
\item [v)] For $k = \frac{5}{2}$, any isomorphism $T:
\mathfrak{D}_{\lambda,\mu}^{\frac{5}{2}} \rightarrow
\mathfrak{D}_{\rho,\nu}^{\frac{5}{2}}$ has a diagonal form by
Proposition \ref{bellie village}. The equivariant conditions of $T$
lead to the following system
\begin{equation}
\label{system5}
\begin{array}{lll}
\tau_0 \beta_0^3(\lambda, \mu) = \tau_3 \beta_0^3(\rho, \nu),\\
\tau_0 \beta_0^4(\lambda, \mu) = \tau_4 \beta_0^4(\rho, \nu),\\
\tau_0 \beta_0^5(\lambda, \mu) = \tau_5 \beta_0^5(\rho, \nu),\\
\tau_1 \beta_1^4(\lambda, \mu) = \tau_4 \beta_1^4(\rho, \nu),\\
\tau_1 \beta_1^5(\lambda, \mu) = \tau_5 \beta_1^5(\rho, \nu),\\
\tau_2 \beta_2^5(\lambda, \mu) = \tau_5 \beta_2^5(\rho, \nu).
\end{array}
\end{equation}
One can readily check that this system has solutions only if
$\lambda = \rho$ or $\rho+ \mu = \frac{1}{2}$. The first isomorphism
is tautological, and the second is just the passage to the adjoint
module.
\item [vi)] For $k > \frac{5}{2}$, 
let $T: \mathfrak{D}_{\lambda,\mu}^k \rightarrow
\mathfrak{D}_{\rho,\nu}^k$ be an isomorphism of
$\mathcal{K}(1)$-modules. The restriction of $T$ to
$\mathfrak{D}_{\lambda,\mu}^{\frac{5}{2}}\subset\mathfrak{D}_{\lambda,\mu}^k$
must be an isomorphism onto $\mathfrak{D}_{\rho,\nu}^{\frac{5}{2}}$.
So, we must have $\lambda = \rho$ or $\rho+ \mu = \frac{1}{2}$.
\end{itemize}
Theorem \ref{main11} is now completely proved for nonresonant values of $\delta=\mu-\lambda$.
In the next subsection, using the approach of the deformation theory
(see, e.g.,\cite{abbo, bbbbk, NS, NIS, nr}), we will give the
relationship between singular modules and cohomology  for
nonresonant values of $\delta$.
\subsection{Cohomological interpretation of singularity of $\mathfrak{D}_{\lambda,\mu}^k$, $k\leq2$}
\label{kinza} Of course, the actions of $\cK(1)$ on ${\mathfrak S}_{\delta}$ defined by $\widetilde{\fL}^{\l,\mu}$ and $\widetilde{\fL}^{\rho,\nu}$ are two $\mathfrak{osp}(1|2)$-trivial deformations of the natural action $\fL$. These deformations become trivial when restricted to $\mathfrak{osp}(1|2)$. So, they are
related to the
the $\mathfrak{osp}(1|2)$-relative cohomology space  \cite{bbbbk}:
\begin{equation*}
\mathrm{H}^1_\mathrm{diff}\left(\mathcal{K}(1), \mathfrak{osp}(1|2);
\mathrm{End}_\mathrm{diff}({\mathfrak
S}_{\delta})\right)=\bigoplus_{p\leq j
}\mathrm{H}^1_\mathrm{diff}\left(\mathcal{K}(1),
\mathfrak{osp}(1|2);
\Pi^{j-p}\left(\mathfrak{D}_{\delta-{j\over2},\delta-{p\over2}}\right)\right),
\end{equation*}
where $\mathrm{H}^1_\mathrm{diff}$ denotes the
differential cohomology; that is, only cochains given by
differential operators are considered. This $\mathfrak{osp}(1|2)$-relative cohomology space is spanned by the nontrivial 1-cocycles $\Pi^{j-p}\left(\pi^{j-p}\circ \frak{J}_{\frac{j-p}{2}+1}^{\delta-\frac{j}{2}}(\cdot~,~)\right)$ where $j-p\in\{3,\,4,\,5,\,6,\,8\}$
(see, e.g., \cite{bbbbk, c}). So, for $k\leq2$, by
fundamental arguments of the theory of deformation \cite{nr}, we can see that, for $j-p=3$ or 4, if $\beta_p^j(\lambda,\mu)=0$ and
$\beta_p^j(\rho,\nu)\neq0$, then
the $\mathcal{K}(1)$-modules $\mathfrak{D}_{\lambda,\mu}^k$ and $\mathfrak{D}_{\rho,\nu}^k$ are not isomorphic. Thus, singular modules appear whenever at least one of the coefficients
$\beta_p^j(\lambda,\mu)$ in (\ref{bellie}) vanishes.

\begin{rmk}
\label{remarque1} Clearly, the $\mathfrak{osp}(1|2)$-trivial deformation of the action
of $\mathcal{K}(1)$ on the space of symbols
$\mathfrak{S}^{\frac{3}{2}}_{\frac{1}{2}-\lambda}$ is trivial. So,
as a $\mathcal{K}(1)$-modules, we have
$\mathfrak{D}_{\lambda,\frac{1}{2}}^{\frac{3}{2}}\simeq\mathfrak{S}^{\frac{3}{2}}_{\frac{1}{2}-\lambda}$.
\end{rmk}
\section{Obstructions to the existence of $\mathfrak{osp}(1|2)$-equivariant symbol mappings}
For the resonant values $\delta$, there exist a series of cohomology
classes of $\mathfrak{osp}(1|2)$ that are obstructions for existence
of the isomorphism in Theorem \ref{main11}. More precisely, consider
the linear mappings
${\Upsilon}_{n}:\mathfrak{osp}(1|2)\rightarrow\mathfrak{D}_{{1-n\over2},{n\over2}}$
given by
\begin{equation}
\label{b3}{\Upsilon}_{n}(X_F)=(-1)^{|F|}\left(
(n-1){\eta}^4(F)\overline{\eta}^{2n-3}
+{\eta}^3(F)\overline{\eta}^{2n-2}\right).
\end{equation}
We can check (see \cite{bb}) that these mappings are nontrivial
1-cocycles on $\mathfrak{osp}(1|2)$ for any
$n\in\mathbb{N}\setminus\{0\}$. Theses cocycles  arises in the
action (\ref{d-action}) of $\mathfrak{osp}(1|2)$ on
${\frak{D}}_{\lambda,\mu}$. We can nevertheless define a canonical
symbol mapping in the resonant case such that its deviation from
$\mathfrak{osp}(1|2)$-equivariance is measured by the corresponding
cocycle (\ref{b3}).
\subsection{$\mathfrak{osp}(1|2)$-modules deformation}
From now on, $\delta\in\{\half, 1, \frac{3}{2}, 2, \dots,k\}$.
Here, we will construct a nontrivial deformation of the natural
action of the Lie superalgebra $\mathfrak{osp}(1|2)$ on ${\frak
S}^k_{\delta}$,
generated by the cocycles (\ref{b3}).
\begin{prop}\label{def act}
The map $\mathcal{L}
: \mathfrak{osp}(1|2)\rightarrow{\rm End}({\frak
S}^k_{\delta})
$
defined by
\begin{equation}
\label{L}
\mathcal{L}_{X_F}\left(\sum_{i=0}^{2k}\Pi^{i}\left(P_i\alpha^{\delta-{i\over2}}\right)\right)
=\sum_{i=0}^{2k}\Pi^{i}\left(\bar{P}^{X_F}_i\alpha^{\delta-{i\over2}}\right)
\end{equation}
with
\[\left\{
\begin{array}{lllllllll}
\bar{P}_i^{X_F}&=&\mathfrak{L}^{\delta-{i\over2}}_{X_F}(P_i)\quad\quad
\hbox{if}
\quad i<4\delta-2k-1\quad\hbox{or}\quad i>2\delta-1,\\[5pt]
\bar{P}_i^{X_F}&=&\mathfrak{L}^{\delta-{i\over2}}_{X_F}(P_i)-\varepsilon_i^s(-1)^{|P_s|}\Big({(s-i-1)\over
2}
{\eta}^4(F)\overline{\eta}^{s-i-2}(P_s)+
{\eta}^3(F)\overline{\eta}^{s-i-1}(P_s)\Big)\\[5pt]
& & \quad\quad\quad\quad\quad\quad\hbox{if}\quad 4\delta-2k-1\leq
i\leq2\delta-1,
\end{array}
\right.
\]
where $s=4\delta-i-1$ and
\begin{equation}
\label{constant} \varepsilon_i^s= \left\{
\begin{array}{lll}
(-1)^{2\delta}\,\left(\lambda+{1\over2}[{s\over2}]\right)\,\gamma^{s-1}_{s-1-i}&\hbox{if}\quad i\quad\hbox{is even} \\[5pt]
-(-1)^{2\delta}\,{s\over4}\, \gamma^{s-1}_{s-1-i}&\hbox{if}\quad
i\quad\hbox{is odd},
\end{array}
\right. \end{equation} is an action of the Lie superalgebra
$\mathfrak{osp}(1|2)$ on the superspace of symbols
${\frak S}^k_{\delta}$ of order $\leq k$.
\end{prop}
\begin{proofname}.
 First, it is easy to see that the map
$
\Gamma:\mathfrak{D}_{\lambda,\mu}\rightarrow\Pi\left(\mathfrak{D}_{\lambda,\mu}\right)
$ defined by $\Gamma(A)=\Pi(\pi\circ A)$ satisfies
\begin{equation*}
\fL^{\lambda,\mu}_{X_F}\circ
\chi=(-1)^{|F|}\Gamma\circ\fL^{\lambda,\mu}_{X_F}\quad\text{ for all
} X_F\in\mathfrak{osp}(1|2).
\end{equation*}
Thus, we deduce the structure of the first cohomology space
$\mathrm{H^1}(\mathfrak{osp}(1|2);\Pi(\mathfrak{D}_{\lambda,\mu}))$
from $\mathrm{H^1}(\mathfrak{osp}(1|2);\mathfrak{D}_{\lambda,\mu})$.
Indeed, to any 1-cocycle $\Upsilon$
on $\mathfrak{osp}(1|2)$ with values in $\mathfrak{D}_{\lambda,\mu}$
corresponds an 1-cocycle $\Gamma\circ\Upsilon$ on
$\mathfrak{osp}(1|2)$ with values in
$\Pi(\mathfrak{D}_{\lambda,\mu})$.
Obviously, $\Upsilon$ is a couboundary if and only if
$\Gamma\circ\Upsilon$ is a couboundary. Second, we can readily check
 that the map $\mathcal{L}$ satisfies the homomorphism condition
\begin{equation*}\label{nizar}
\mathcal{L}_{[X_F,X_G]}=[\mathcal{L}_{X_F},
\mathcal{L}_{X_G}]\quad\text{ for all }
X_F,\,X_G\in\mathfrak{osp}(1|2).
\end{equation*}
So, the map $\mathcal{L}$ is the nontrivial deformation of the
natural action of 
$\mathfrak{osp}(1|2)$ on ${\frak S}^k_{\delta}$ generated by the
cocycles (\ref{b3}), up to the map $\Gamma$.\hfill$\Box$
\end{proofname}

\medskip

Denote by $\mathcal{M}^k_{\lambda, \mu}$ the
$\mathfrak{osp}(1|2)$-module structure on ${\frak
S}^k_{\delta}$ defined by $\mathcal{L}$ for a fixed $\lambda$ and
$\mu$.
\begin{rmk} Note that the map $\mathcal{L}$ given in Proposition \ref{def act} define an action of $\mathfrak{osp}(1|2)$ on ${\frak
S}^k_{\delta}$ for any scalars replacing those in (\ref{constant}).
\end{rmk}
\subsection{Normal symbol} 
Here, we prove existence and uniqueness (up to normalization) of
$\mathfrak{osp}(1|2)$-isomorphism between
${\frak{D}}_{\lambda,\mu}^k$ and $\mathcal{M}^k_{\lambda, \mu}$
providing a $\lq\lq$ total symbol" of differential operators in the
resonant cases. The following Proposition gives the existence of such
an isomorphism.
\begin{prop}
\label{act} 
There exists an $\mathfrak{osp}(1|2)$-invariant symbol map called a
normal symbol map
\begin{equation}
\label{jamila} \widetilde{\sigma}_{\lambda,
\mu}:{\frak{D}}_{\lambda,\mu}^k\stackrel{\simeq}\longrightarrow\mathcal{M}^k_{\lambda,
\mu}.
\end{equation}
It sends a differential operator $A=\sum_{i=0}^{2k} a_i(x,\theta)\,
\overline{\eta}^i$ to the tensor density
\begin{equation}
\label{jamila1} \widetilde{\sigma}_{\lambda,
\mu}(A)=\sum^{2k}_{j=0}\Pi^{j}\left(\widetilde{a}_j\alpha^{\delta-{j\over
2}}\right),
\end{equation}
where $\widetilde{a}_j=\sum^{2k}_{i\geq j}\xi^i_j\,\eta^{i-j}(a_i)$
with
\begin{equation}
\left\{
\begin{array}{lllllllll}
\xi^i_j&=&\omega^{i}_{j,s}\xi^s_j +\kappa_{j,s}^i&\hbox{if}\quad 4\delta-2k-1\leq j\leq2\delta-1\hbox{ and } i>s,\\[5pt]
\xi^i_j&=&\gamma^i_{i-j} &\hbox{otherwise},
\end{array}
\right.
\label{system1}
\end{equation}
where 
\begin{equation}
\label{benfraj} \omega^{i}_{j,s}=(-1)^{\lbrack
\frac{i-s}{2}\rbrack}\frac{\left(\begin{array}{c}
\lbrack \frac{i}{2}\rbrack\\
\lbrack \frac{2(i-s)+1-(-1)^{(i+1)s}}{4}\rbrack \end{array}\right)
\left(\begin{array}{c}\lbrack \frac{i-1}{2}\rbrack+2\lambda\\
\lbrack \frac{2(i-s)+1-(-1)^{i(s+1)}}{4}\rbrack
\end{array}\right)}{ \left(\begin{array}{c}
\lbrack \frac{i-j}{2}\rbrack\\
\lbrack \frac{2(i-s)+1+(-1)^{i+j}}{4}\rbrack \end{array}\right)
\left(\begin{array}{c}2\delta-\lbrack\frac{j+s}{2}\rbrack-2\\
\lbrack \frac{i-s}{2}\rbrack \end{array}\right)}
\end{equation}
and
\[
\kappa_{j,s}^i=(-1)^{2\delta}\varepsilon_j^s\sum_{\ell=s+1}^i\frac{\omega^{i}_{j,\ell}}{\vartheta_{j,\ell}}\gamma^\ell_{\ell-s}
\]
with $\gamma^i_{i-j}$ is as (\ref{valentain}), $s=4\delta-j-1$, and
\[
\vartheta_{j,\ell}=\left\{
\begin{array}{llllll}
\frac{1}{2}[\frac{\ell-j}{2}]+(\delta-\frac{\ell}{2})&\hbox{if}~~\ell-j~~\hbox{is
odd}\\[5pt]
\frac{\ell-j}{4}&\hbox{if}~~\ell-j~~\hbox{is even}.
\end{array}
\right.
\]
\end{prop}

\begin{proofname}. 
For $A\in {\frak{D}}_{\lambda,\mu}^k$ and
$X_F\in\mathfrak{osp}(1|2)$, we have {\small\begin{equation}
\label{nisma}
 \widetilde{\sigma}_{\lambda,
\mu}\left(\fL^{\l,\mu}_{X_F}(A)\right)=\sum^{2k}_{j=0}\Pi^{j}\left(\widetilde{a}_j^{X_F}\alpha^{\delta-{j\over
2}}\right)\quad\text{with}\quad\widetilde{a}^{X_F}_j=\sum^{2k}_{i\geq
j}\xi^i_j\eta^{i-j}(a^{X_F}_i).
\end{equation}}
Substituting expression (\ref{naturel action}) for $a^{X_F}_i$
in (\ref{nisma}), we get
\begin{equation}
\renewcommand{\arraystretch}{1.4}
\begin{array}{llll}
\widetilde{a}^{X_F}_j&=&\frak{L}^{\delta-{j\over2}}_{X_F}(\widetilde{a}_j)- \displaystyle\sum_{i\geq
j+1}^{2k}\rho^i_j\eta(F')\overline{\eta}^{i-j-1}(a_i)
-\displaystyle\sum_{i\geq
j+2}^{2k}\widetilde{\rho}^i_jF''\overline{\eta}^{i-j-2}(a_i),
\end{array}
\label{proof_equation1}
\end{equation}
where
{\small
\[
\begin{array}{llll}
\rho^i_j=&(-1)^{(i-j)(|a_i|+1)}\left(\Lambda^i_j\zeta_{i-j-1,i-j,\delta-{i\over2}}\xi^i_j
+
\Lambda^i_{j+1}\zeta_{i-1,i,\lambda}\xi^{i-1}_j\right),\\
\widetilde{\rho}^i_j=&(-1)^{i-j)(|a_i|+1)}\Big(\Lambda^i_j\zeta_{i-j-2,i-j,\delta-{i\over2}}\xi^i_j-
{\small{\begin{pmatrix}i-j-1\\1\end{pmatrix}}_s}
\Lambda^i_{j+1}\zeta_{i-1,i,\lambda}\xi^{i-1}_j-\\& \Lambda^i_{j+2}\zeta_{i-2,i,\lambda}\xi^{i-2}_j\Big)
\end{array}
\]}
with $\Lambda^i_j=(-1)^{[{i-j\over2}]}$ and $\zeta_{i,j,\lambda}$ is as in (\ref{zeta}).
So, we can see that, for $j<4\delta-2k-1$ or $ j>2\delta-1$, the
symbol map (\ref{jamila}) commutes with the action of
$\mathfrak{osp}(1|2)$ if and only if the following system is
satisfied: {\small\begin{equation} \label{benfraj1}
\begin{array}{llll}
{(i-j)\over2}\xi^i_j={i\over2}\xi^{i-1}_j&\hbox{if}~i~\hbox{and}~j~\hbox{are
even},\\[10pt]
{(i-j)\over2}\xi^i_j=\left(2\lambda+\lbrack{i\over2}\rbrack\right)\xi^{i-1}_j&\hbox{if}~i~\hbox{and}~j~\hbox{are odd},\\[10pt]
\left(i-2\delta-\lbrack{i-j\over2}\rbrack\right)\xi^i_j=\left(2\lambda+\lbrack{i\over2}\rbrack\right)
\xi^{i-1}_j&\hbox{if}~i~\hbox{is odd and}~j~\hbox{is even},\\[10pt]
\left(i-2\delta-\lbrack{i-j\over2}\rbrack\right)\xi^i_j={i\over2}
\xi^{i-1}_j&\hbox{if}~i~\hbox{is even and}~j~\hbox{is odd}.
\end{array}
\end{equation}}
Hence, the solution of the system (\ref{benfraj1}) with the initial
condition $\xi^j_j=1$ is unique
and given by $\xi^i_j=\gamma^i_{i-j}$. 
For $4\delta-2k-1\leq j\leq2\delta-1$, the
$\mathfrak{osp}(1|2)$-equivariance condition reads 
{\small\begin{equation}
\widetilde{a}^{X_F}_j
=\frak{L}^{\delta-{j\over2}}_{X_F}(\widetilde{a}_j)-\varepsilon_j^s(-1)^{|\widetilde{a}_s|}\Big({(s-j-1)\over
2}
{\eta}^4(F)\overline{\eta}^{s-j-2}(\widetilde{a}_s)+
{\eta}^3(F)\overline{\eta}^{s-j-1}(\widetilde{a}_s)\Big).
\label{proof_equation2}
\end{equation}}
Thus, it is easy to see that the 
solutions of equation (\ref{proof_equation2}) with indeterminate
$\xi_j^i$  are given by (\ref{system1}).\hfill$\Box$
\end{proofname}

\medskip

To study the uniqueness (up to normalization) for the symbol map
given by (\ref{jamila}), we need the following result.
\begin{prop}
\label{mabrouk} The action of $\mathcal{K}(1)$ over
$\mathcal{M}^k_{\lambda, \mu}$ in terms of
$\mathfrak{osp}(1|2)$-equivariant normal symbols is given by:
{\small\begin{equation} \label{mabrouk ben ammar}
\begin{array}{llllllllllll}
\widetilde{a}_p^{X_F} &=& \fL^{\delta - \frac{p}{2}}_{X_F}
(\widetilde{a}_p) +
\sum_{j=p+3}^{2k}\chi_p^j\,
\pi^{j-p}\circ\frak{J}_{\frac{j-p}{2}+1}^{\delta - \frac{j}{2}
}\left(X_F,\widetilde{a}_j\right)& \hbox{ if }
p<4\delta-2k-1\\[5pt]& & ~& \hbox{ or } p>2\delta-1,\\[10pt]
\widetilde{a}_p^{X_F}&=&\mathfrak{L}^{\delta-{p\over2}}_{X_F}(\widetilde{a}_p)+\sum_{j=p+3}^{2k}\chi_p^j\,
\pi^{j-p}\circ\frak{J}_{\frac{j-p}{2}+1}^{\delta - \frac{j}{2}
}\left(X_F,\widetilde{a}_j\right)\\[8pt]
&-&\varepsilon_p^s(-1)^{|\widetilde{a}_s|}\Big({(s-p-1)\over 2}
{\eta}^4(F)\overline{\eta}^{s-p-2}(\widetilde{a}_s)+
{\eta}^3(F)\overline{\eta}^{s-p-1}(\widetilde{a}_s)\Big)\\&
&~&\hbox{if}\quad 4\delta-2k-1\leq p\leq2\delta-1\\[5pt]& & ~& \hbox{ and } s\neq p+5,\\[10pt]
\widetilde{a}_p^{X_F}&=&\mathfrak{L}^{\delta-{p\over2}}_{X_F}(\widetilde{a}_p)+\sum_{j=p+3}^{2k}\chi_p^j\,
\pi^{j-p}\circ\frak{J}_{\frac{j-p}{2}+1}^{\delta - \frac{j}{2}
}\left(X_F,\widetilde{a}_j\right)\\[8pt]
&-&\varepsilon_p^s(-1)^{|\widetilde{a}_s|}\Big({(s-p-1)\over 2}
{\eta}^4(F)\overline{\eta}^{s-p-2}(\widetilde{a}_s)+
{\eta}^3(F)\overline{\eta}^{s-p-1}(\widetilde{a}_s)\Big)\\[8pt]
&+&\Xi_p^s(-1)^{|\widetilde{a}_s|}\Big(4
{\eta}(F''')\widetilde{a}_s+
F'''\overline{\eta}(\widetilde{a}_s)\Big)\\&
&~&\hbox{if}\quad 4\delta-2k-1\leq p\leq2\delta-1\\[5pt]& & ~& \hbox{ and } s= p+5,
\end{array}
\end{equation}}
where $\chi_p^j $ and $\Xi_p^s$ are functions of $(\lambda,\mu)$ and
$\varepsilon_p^s$ is as in (\ref{constant}) with  $s=4\delta-p-1$.
\end{prop}
To prove Proposition \ref{mabrouk}, we need the following classical fact:
\begin{lem}
\label{benfra}
Consider a linear differential operator
$b:\mathcal{K}(1)\longrightarrow\frak{D}_{\lambda,\mu}$.

(i) If 
$b$ is an $1$-cocyle
vanishing on $\frak{osp}(1|2)$, then  $b$ is a supertransvectant.

(ii) If $b$
satisfies
\begin{equation*}
\Delta(b)(X,Y)=b(X)=0\hbox{ for all }{X\in\frak {osp}}(1|2),
\end{equation*}
where $\Delta$ stands for differential of cochains on
$\mathcal{K}(1)$ with values in $\frak{D}_{\lambda,\mu}$ (see, e.g.,
\cite{bbbbk, bbo, c}), then $b$ is a supertransvectant.
\end{lem}
We also need the following
\begin{theorem}\cite{bbbbk,c}
\label{cbbbbk}
\label{th1} ${\rm dim}\mathrm{H}^1_{\rm
diff}(\mathcal{K}(1);\mathfrak{D}_{\lambda,\mu})=1$ if
\begin{equation*}
\begin{array}{llll}
\mu-\lambda=0 &\hbox{ for all } \lambda,\\
\mu-\lambda=\frac{3}{2}&\hbox{ for all }~\lambda,\\ \mu-\lambda=2
&\hbox{ for all } \lambda,\\ \mu-\lambda=\frac{5}{2}&\hbox{ for
all }~\lambda,\\ \mu-\lambda=3 &\hbox{ and }
\lambda\in\{0,\,-\frac{5}{2}\},\\ \mu-\lambda=4 &\hbox{ and }
\lambda=\frac{-7\pm\sqrt{33}}{4}.
\end{array}
\end{equation*}
${\rm dim}\mathrm{H}^1_{\rm
diff}(\mathcal{K}(1);\frak{D}_{0,\frac{1}{2}})=2$. Otherwise,
$\mathrm{H}^1_{\rm
diff}(\mathcal{K}(1);\frak{D}_{\lambda,\mu})=0$.

The spaces $\mathrm{H}^1_{\rm
diff}(\mathcal{K}(1);\frak{D}_{\lambda,\mu})$ are spanned by the
cohomology classes of of
$\Upsilon_{\lambda,\lambda+\frac{k}{2}}=\frak{J}_{\frac{k}{2}+1}^{\lambda}$,
where $k\in\{3,\,4,\,5,\,6,\,8\}$, and by the cohomology classes of the
following $1$-cocycles:
\begin{equation*}\left\{\begin{array}{llllll}
\Upsilon_{\lambda,\lambda}(X_F)(G\alpha^{\lambda})&=&F'G\alpha^{\lambda},\\[4pt]
\Upsilon_{0,\frac{1}{2}}(X_F)(G)&=&
\overline{\eta}(F') G\alpha^{{1\over2}},\\[4pt]
\widetilde{\Upsilon}_{0,\frac{1}{2}}(X_F)( G)&=&
\overline{\eta}(F' G)\alpha^{{1\over2}},\\[4pt]
\Upsilon_{-\frac{1}{2},1}(X_F)( G\alpha^{-{1\over2}})&=&\Big(
\overline{\eta}(F')G'+(-1)^{|F|}F''\overline{\eta}(G)\Big)\alpha\\[4pt]
\Upsilon_{-1,\frac{3}{2}}(X_F)( G\alpha^{-1})&=&
\Big((-1)^{|F|}(F'''\overline{\eta}(G)+2F''\overline{\eta}(G'))+
2\overline{\eta}(F'')G'+\overline{\eta}(F')G''\Big)\alpha^{\frac{3}{2}}.
\end{array}\right.\end{equation*}
\end{theorem}
\begin{proofname} (Proposition \ref{mabrouk}). By direct computation, using formula (\ref{naturel action}),
we can see that the action 
of $\mathcal{K}(1)$ over
$\mathcal{M}^k_{\lambda, \mu}$ in terms of
$\mathfrak{osp}(1|2)$-equivariant normal symbols is given by:
\[
\widetilde{a}_p^{X_F} = \fL^{\delta - \frac{p}{2}}_{X_F} (\widetilde{a}_p) + (\hbox{terms
in}~~ \overline{\eta}^{n}(F),~n\geq3).
\]
So, it is a deformation of the natural action of $\mathcal{K}(1)$
over $\mathfrak{S}^k_{\delta}$ such that its restriction to
$\mathfrak{osp}(1|2)$ coincides with the map $\mathcal{L}$ given by
(\ref{L}). Therefore, according to Theorem \ref{cbbbbk} and Lemma
\ref{benfra}, we deduce that
it is given by:

\medskip

\noindent For $p<4\delta-2k-1$ or $ p>2\delta-1$,
{\small\begin{equation}
\label{n1}
\widetilde{a}_p^{X_F} = \fL^{\delta - \frac{p}{2}}_{X_F}
(\widetilde{a}_p) + \sum_{j=p+3}^{2k}\hbar_p^j\,
\pi^{j-p}\circ\frak{J}_{\frac{j-p}{2}+1}^{\delta - \frac{j}{2}
}\left(X_F,\widetilde{a}_j\right)+\sum_{j=0}^{2k}C_{j,p}\left(X_F,\widetilde{a}_j\right).
\end{equation}}

\medskip

\noindent For $4\delta-2k-1\leq p\leq2\delta-1$  and  $s\neq p+5$,
{\small\begin{equation}
\label{n2}
\begin{array}{lllllll}
\widetilde{a}_p^{X_F}&=&\mathfrak{L}^{\delta-{p\over2}}_{X_F}(\widetilde{a}_p)+\sum_{j=p+3}^{2k}\hbar_p^j\,
\pi^{j-p}\circ\frak{J}_{\frac{j-p}{2}+1}^{\delta - \frac{j}{2}
}\left(X_F,\widetilde{a}_j\right)+\sum_{j=0}^{2k}C_{j,p}\left(X_F,\widetilde{a}_j\right)\\[8pt]
&-&\varepsilon_p^s(-1)^{|\widetilde{a}_s|}\Big({(s-p-1)\over 2}
{\eta}^4(F)\overline{\eta}^{s-p-2}(\widetilde{a}_s)+
{\eta}^3(F)\overline{\eta}^{s-p-1}(\widetilde{a}_s)\Big).
\end{array}
\end{equation}}

\medskip

\noindent For $4\delta-2k-1\leq p\leq2\delta-1$  and  $s= p+5$,
{\small\begin{equation}
\label{n3}
\begin{array}{lllllll}
\widetilde{a}_p^{X_F}&=&\mathfrak{L}^{\delta-{p\over2}}_{X_F}(\widetilde{a}_p)+\sum_{j=p+3}^{2k}\hbar_p^j\,
\pi^{j-p}\circ\frak{J}_{\frac{j-p}{2}+1}^{\delta - \frac{j}{2}
}\left(X_F,\widetilde{a}_j\right)+\sum_{j=0}^{2k}C_{j,p}\left(X_F,\widetilde{a}_j\right)\\[8pt]
&+&\Re_p^s(-1)^{|\widetilde{a}_s|}\Big(2
F''\overline{\eta}^{3}(\widetilde{a}_s)+
{\eta}^3(F)\overline{\eta}^{4}(\widetilde{a}_s)-F'''\overline{\eta}(\widetilde{a}_s)-
2{\eta}(F'')\overline{\eta}^{2}(\widetilde{a}_s)\Big)\\[8pt]
&+&\aleph_p^s(-1)^{|\widetilde{a}_s|}\Big(2
F''\overline{\eta}^{3}(\widetilde{a}_s)+
{\eta}^3(F)\overline{\eta}^{4}(\widetilde{a}_s)\Big),
\end{array}
\end{equation}}
\noindent where $\hbar_p^j, \Re_p^s, \aleph_p^s$ are functions of $(\lambda,\mu)$ and $C_{j,p}:\mathcal{K}(1)\rightarrow\frak{D}_{\delta-{j\over2},\delta-{p\over2}}$ are linear maps vanishing on $\mathfrak{osp}(1|2)$ with the same parity as the integer
$j-p$. The homomorphism condition of the action of $\mathcal{K}(1)$ over
$\mathcal{M}^k_{\lambda, \mu}$ in terms of
$\mathfrak{osp}(1|2)$-equivariant normal symbols
implies that $\Delta C_{j,p}$ vanish on $\mathfrak{osp}(1|2)$. So, by Lemma \ref{benfra}, we can see that the maps $C_{j,p}$ are supertransvectants vanishing on  $\mathfrak{osp}(1|2)$. Therefore, we deduce that formulas (\ref{n1}) and (\ref{n2}) can be expressed as in Proposition \ref{mabrouk}. Moreover, for $s=p+5$, up to a scalair factor, the supertransvectant $\frak{J}_{\frac{7}{2}}^{-1}$ is given by
\begin{equation}
\label{n4}
\frak{J}_{\frac{7}{2}}^{-1}\left(X_F,\widetilde{a}_s\right)=(-1)^{|F|}\Big(2
{\eta}(F''')(\widetilde{a}_s)-3
{\eta}(F'')\overline{\eta}^{2}(\widetilde{a}_s)-F'''\overline{\eta}(\widetilde{a}_s)\Big).
\end{equation}
Using formula (\ref{n4}), we deduce that expression (\ref{n3}) can
be expressed as in Proposition \ref{mabrouk}. This completes the
proof.\hfill$\Box$
\end{proofname}


\begin{cor}
\label{ben ammar} The constants $\chi_p^j$ given by (\ref{mabrouk
ben ammar}) and $\xi^j_p$ given by (\ref{system1}) satisfy the
following relations
\begin{equation} \label{new}
\begin{array}{llllllllllll}
\Theta^j_p\chi_p^j&=&\varsigma^j_p-\sum_{r=p+3}^{j-1}(-1)^{(j-p)(j-r)}\Lambda^j_r\Theta^r_p\chi_p^r\xi^j_r\quad\hbox{if}~j\geq p+4,\\[5pt]
\Theta^{p+3}_p\chi_p^{p+3}&=& \varsigma^{p+3}_p,
\end{array}
\end{equation}
where \begin{equation} \label{new1}\begin{array}{llllllllllll}
(-1)^{j-p}\varsigma^j_p&=&-\Lambda^j_p\zeta_{j-p-3,j-p,\delta-{j\over2}}\xi_p^j
-\Lambda^{j-1}_p{\begin{pmatrix}j-p-1\\2\end{pmatrix}}_s\zeta_{j-1,j,\lambda}\xi_p^{j-1}\\[7pt]& &+
\Lambda^{j-2}_p{\begin{pmatrix}j-p-2\\1\end{pmatrix}}_s\zeta_{j-2,j,\lambda}\xi_p^{j-2}
-\Lambda^{j-3}_p\zeta_{j-3,j,\lambda}\xi_p^{j-3},
\end{array}
\end{equation}
$\zeta_{i,j,\lambda}$ is as in (\ref{zeta}) and $\Theta^j_p$
is the coefficient of
$$(-1)^{|F|+(j-p)|\widetilde{a}_j|}\,\overline{\eta}(F'')\overline{\eta}^{j-p-3}(\widetilde{a}_j)\quad
\hbox{in}~~ \pi^{j-p}\circ\frak{J}_{\frac{j-p}{2}+1}^{\delta -
\frac{j}{2} }\left(X_F,\widetilde{a}_j\right).$$
\end{cor}
\begin{proofname}. 
First, using the fact that
$\widetilde{a}_r=\sum^{2k}_{i=r}\xi^i_r\eta^{i-r}(a_i)$ 
and formula (\ref{mabrouk ben ammar}), we can check that the
coefficient of
$(-1)^{|F|+(j-p)|a_j|}\overline{\eta}(F'')\overline{\eta}^{j-p-3}({a}_j)
~~\hbox{ in }~ \widetilde{a}^{X_F}_p~\hbox{ for }~ j\geq p+3,$ is
\[
\Theta^j_p\chi_p^j+\sum_{r=p+3}^{j-1}(-1)^{(j-p)(j-r)}\Lambda^j_r\Theta^r_p\chi_p^r\xi^j_r.
\]
On the other hand, we have
$\widetilde{a}^{X_F}_p=\sum^{2k}_{i=p}\xi^i_p\eta^{i-p}(a^{X_F}_i)$,
where $a^{X_F}_i$ is given by (\ref{naturel action}). Thus, by
direct computation, we can see that the coefficient of
$(-1)^{|F|+(j-p)|a_j|}\overline{\eta}(F'')\overline{\eta}^{j-p-3}({a}_j)$
in $\widetilde{a}^{X_F}_p$ is $\varsigma^j_p$ given by (\ref{new1}).
Corollary \ref{ben ammar} is proved.\hfill$\Box$
\end{proofname}

\medskip

The normal symbol map depends on the choice of $\xi^s_p$, which play
a role arbitrary. Clearly, $s-p$ is odd. Moreover, we can readily
check that the coefficient of $\xi^s_p$ in $\chi_p^{s+2}$ vanishes
for $s=p+1$. So, in the following, we will use the normal symbol map
uniquely defined, up to a scalar factor, by imposing the following
condition to $\xi^s_p$:
\begin{itemize}
\item[(i)] if $s\geq p+3$ we choose $\xi^s_p$ such that $\chi^s_p=0$,

\item[(ii)] if $s= p+1$ we choose $\xi^s_p$ so as to cancel the first
term of the following sequence, where the coefficient of $\xi^s_p$
is nonzero:
\begin{equation}
\label{condition}
\chi_p^{s+3},~\chi_{p-3}^s,~\chi_{p-4}^s,\ldots,\chi_{0}^s.
\end{equation}
\end{itemize}
Note that this choice is possible thanks to Corollary \ref{ben
ammar}.

\section{Proof of Theorem \ref{main11} in the resonant case}
\label{proof2}
 The existence and uniqueness of the normal symbol allow us, by a similar
process to that used in section \ref{proof1}, to complete the proof
of Theorem \ref{main11}.
\begin{prop}
\label{bellie village1} Let $T: \mathfrak{D}_{\lambda,\mu}^k
\rightarrow \mathfrak{D}_{\rho,\nu}^k$
 be an isomorphism of $\mathcal{K}(1)$-modules. Then $T$ is diagonal
 in terms of normal symbols.
 \end{prop}
\begin{proofname}. Similar to that of Proposition \ref{bellie
village}.\hfill$\Box$
\end{proofname}

Now, let $A\in\mathfrak{D}_{\lambda,\mu}^k$. The normal symbol of
$T(A)$ is
\[
\widetilde{\sigma}_{\lambda,
\mu}(T(A))=\sum^{2k}_{j=0}\Pi^{j}\left(\widetilde{a}^T_j\alpha^{\delta-{j\over
2}}\right).
\]
Proposition \ref{bellie village1} implies that there exist a
constants $\tau_0,\ldots, \tau_{2k}$ depending on $\lambda, \mu,
\rho$ and $\nu$, such that $\widetilde{a}^T_j=\tau_j\widetilde{a}_j$
for all $j=0,\ldots, 2k$. The condition of
$\mathfrak{osp}(1|2)$-equivariance of $T$   in terms of normal
symbol, leads to the following system:
\begin{equation} \label{new3}
\begin{array}{llllll}
&\tau_p\,\chi_p^{j}(\lambda, \mu)=\tau_j\,\chi_p^{j}(\rho,
\nu)&\hbox{for}\quad p=0,\ldots, 2k~~\hbox{and}~j\geq
p+3,\\[5pt]
&\tau_p\,\varepsilon_p^s(\lambda, \mu)=\tau_s\,\varepsilon_p^s(\rho,
\nu)&\hbox{for}\quad 4\delta-2k-1\leq
p\leq2\delta-1,\\[5pt]
&\tau_p\,\Xi_p^s(\lambda, \mu)=\tau_s\,\Xi_p^s(\rho,
\nu)&\hbox{for}\quad 4\delta-2k-1\leq
p\leq2\delta-1~~\hbox{and}~s=p+5,
\end{array}
\end{equation}
where the $\chi_p^{j}$ are given by (\ref{mabrouk ben ammar}) for
$j\geq p+3$, and $\varepsilon_p^s$ is as in (\ref{constant}) with
$s=4\delta-p-1$.
\subsection{Isomorphisms of $\mathcal{K}(1)$-modules in terms of normal
symbol}
The resolution of the system  (\ref{new3}) shows that the
isomorphisms of  $\mathcal{K}(1)$-modules in terms of normal symbol,
in the resonant case, are an extension, except for $(k, \delta)=(2,
2)$, of isomorphisms in terms of $\mathfrak{osp}(1|2)$-equivariant
symbols in the nonresonant case. Indeed, by solving the system
(\ref{new3}), using formulas (\ref{new}) and (\ref{system1}) with
the help of condition (\ref{condition}), we get:

\begin{itemize}
  \item [i)]
 For $k = \frac{1}{2}$, an
isomorphism $T: \mathfrak{D}_{\lambda,\mu}^{\frac{1}{2}} \rightarrow
\mathfrak{D}_{\rho,\nu}^{\frac{1}{2}}$
is obtained by taking
\begin{equation*}
 \left(\widetilde{a}_0^T, \widetilde{a}_1^T\right) =
\left(\frac{\rho}{\lambda} \,\widetilde{a}_0, \widetilde{a}_1\right).
\end{equation*}
\item [ii)] For $k = 1$, an
isomorphism $T: \mathfrak{D}_{\lambda,\mu}^1 \rightarrow
\mathfrak{D}_{\rho,\nu}^1$
is obtained by taking (with $\tau\neq0$)
\[
\left\{
\begin{array}{llllllll}
\left(\widetilde{a}_0^T, \widetilde{a}_1^T, \widetilde{a}_2^T\right) =
\left(\frac{\rho}{\lambda}\, \widetilde{a}_0,~ \widetilde{a}_1,~ \tau
\widetilde{a}_2\right)&\hbox{for} ~~\delta=\frac{1}{2},\\[5pt]
\left(\widetilde{a}_0^T, \widetilde{a}_1^T, \widetilde{a}_2^T\right) =
(\tau\widetilde{a}_0, \widetilde{a}_1, \widetilde{a}_2)&\hbox{for}
~~\delta=1.
\end{array}
\right.
\]
\item [iii)] For $k = \frac{3}{2}$, we get
a
an isomorphism $T: \mathfrak{D}_{\lambda,\mu}^{\frac{3}{2}}
\rightarrow \mathfrak{D}_{\rho,\nu}^{\frac{3}{2}}$ by taking in
(\ref{new3}):
\[
\begin{array}{lllllllllll}
\tau_0= \frac{\rho^2}{\lambda^2},& \tau_1 =
\frac{\rho}{\lambda},& \tau_3 = 1,& \tau_2\neq0&\hbox{for} ~~\delta=\frac{1}{2},\\[5pt]
\tau_0= \frac{\rho(2\rho+1)}{\lambda(2\lambda+1)},& \tau_1=\tau_2\neq0,& \tau_3=1,& &\hbox{for} ~~\delta=1,\\[5pt]
\tau_0= \frac{\rho(\rho+1)}{\lambda(\lambda+1)},&\tau_2=
\frac{2\rho+1}{2\lambda+1},& \tau_3=1,& \tau_1\neq0&\hbox{for}
~~\delta=\frac{3}{2}.
\end{array}
\]
\item [iv)] For $k = 2$, we get
a
an isomorphism $T: \mathfrak{D}_{\lambda,\mu}^2 \rightarrow
\mathfrak{D}_{\rho,\nu}^2$ by taking in (\ref{new3}):
\[
\begin{array}{lllllllllll}
\tau_0= \frac{\rho^2}{\lambda^2},& \tau_1 =
\frac{\rho}{\lambda},& \tau_3 =\tau_4= 1,& \tau_2\neq0&\hbox{for} ~~\delta=\frac{1}{2},\\[5pt]
\tau_0=
\frac{2\rho+1}{2\lambda+1},&\tau_1=\tau_2=\frac{4\rho+1}{4\lambda+1},&\tau_3=\frac{\lambda}{\rho},&
\tau_4=1&\hbox{for} ~~\delta=1,\\[5pt]
\tau_0= \frac{\rho(\rho+1)}{\lambda(\lambda+1)},&\tau_2=
\frac{2\rho+1}{2\lambda+1},& \tau_1=\tau_3=\tau_4=1,& &\hbox{for}
~~\delta=\frac{3}{2}.
\end{array}
\]
The case $\delta=2$ is particularly because the isomorphisms of
nonresonant case do not extend to the resonant case. Indeed, the
equivariance condition of an isomorphism $T:
\mathfrak{D}_{\lambda,\mu}^2 \rightarrow \mathfrak{D}_{\rho,\nu}^2$
implies
\[
\begin{array}{lllllllllll}
\frac{\tau_0}{\tau_3}=\frac{\chi^3_0(\rho, \nu)}{\chi^3_0(\lambda,
\mu)}= \frac{\rho(2\rho+3)}{\lambda(2\lambda+3)},\quad
\frac{\tau_1}{\tau_4}=\frac{\chi^4_0(\rho, \nu)}{\chi^4_0(\lambda,
\mu)}= \frac{4\rho+3}{4\lambda+3},\quad\frac{\tau_3}{\tau_4}=
\frac{\chi^4_3(\rho, \nu)}{\chi^4_3(\lambda, \mu)}=1,\\[5pt]
\frac{\tau_0}{\tau_4}=\frac{\chi^4_0(\rho, \nu)}{\chi^4_0(\lambda,
\mu)}=
\frac{\rho(2\rho+3)(2\rho^2+3\rho+1)(\lambda^2+\frac{3}{2}\lambda+1)}
{\lambda(2\lambda+3)(2\lambda^2+3\lambda+1)(\rho^2+\frac{3}{2}\rho+1)}.
\end{array}
\]
This system has a solution if $\lambda=\rho$ or
$\rho+\lambda=\frac{1}{2}$.
\item [v)] For $k = \frac{5}{2}$, the system (\ref{new3}) has a solution if $\lambda=\rho$ or
$\rho+\lambda=\frac{1}{2}$. The first isomorphism is tautological,
and the second is just the passage to the adjoint module. The case
$k > \frac{5}{2}$ is deduced from the case $k = \frac{5}{2}$.
\end{itemize}
Now, 
for fixed $k\leq 2$, the same arguments as in subsection \ref{kinza}
together with Theorem \ref{cbbbbk} show that singular modules appear
whenever at least one of the coefficients $\chi_p^j$ or
$\varepsilon_p^s$ in (\ref{mabrouk ben ammar}) vanishes. Theorem
\ref{main11} is proved for resonant case.\hfill$\Box$
\begin{rmk}
For $k=2$, the resonant case $\delta=2$ seems to be particularly
interesting. 
\end{rmk}

\section{Differential Operators on ${S}^{1|n}$}
In this section, we consider the supercircle ${S}^{1|n}$ instead of
${S}^{1|1}$. That is, we consider the supercircle ${S}^{1|n}$ for
$n\geq2$ with local coordinates $(x,\theta_1,\dots,\theta_n),$ where
$\theta=(\theta_1,\dots,\theta_n)$ are odd variables. Any contact
structure on ${S}^{1|n}$ can be reduced to a canonical one, given by
the following  $1$-form:
\begin{equation*}
\a_n=dx+\sum_{i=1}^n\theta_i d\theta_i.
\end{equation*}
The space of $\lambda$-densities will be denoted
\begin{equation}
\label{densities} \mathfrak{F}_\l^n=\left\{F(x,\theta)\a_n^{\l} \mid
F(x,\theta) \in C^\infty({S}^{1|n})\right\}.
\end{equation}
We denote by $\mathfrak{D}^n_{\lambda,\mu}$ the space of
differential operators from $\mathfrak{F}_\lambda^n$ to
$\mathfrak{F}_\mu^n$ for any $\lambda,\,\mu\in\mathbb{R}$. The Lie
superalgebra $\mathcal{K}(n)$ of contact vector fields on
${S}^{1|n}$ is spanned by the vector fields of the form (see,e.g;
\cite{NBS, gls}):
\begin{equation*}
X_F=F\partial_x
-\frac{1}{2}(-1)^{|F|}\sum_i\overline{\eta}_i(F)\overline{\eta}_i,\quad\text{where}\quad
F\in C^\infty({S}^{1|n}).
\end{equation*}
where $\overline{\eta}_i=\partial_{\theta_i}-\theta_i\partial_x$.
Since $-\eta_i^2=\partial_x,$ and
$\partial_i=\eta_i-\theta_i\eta_i^2,$ every differential operator
$A\in\mathfrak{D}^n_{\lambda,\mu}$ can be expressed in the form
\begin{equation}
\label{diff1}
A(F\alpha_n^\lambda)=\sum_{\ell=(\ell_1,\dots,\ell_n)}a_\ell(x,\theta)
\eta_1^{\ell_1}\dots\eta_n^{\ell_n}(F)\alpha_n^\mu,
\end{equation}
where the coefficients $a_\ell(x,\theta)\in C^\infty({S}^{1|n})$
(see \cite{bbo}). For $k\in\frac{1}{2}\mathbb{N}$, we denote by
$\mathfrak{D}^{n,k}_{\lambda,\mu}$ the subspace of
$\mathfrak{D}^n_{\lambda,\mu}$ of the form
\begin{equation}
\label{diff2}
A(F\alpha_n^\lambda)=\sum_{\ell_1+\cdots+\ell_n\leq2k}a_{\ell_1,\dots,\ell_n}(x,\theta)
\eta_1^{\ell_1}\dots\eta_n^{\ell_n}(F)\alpha_n^\mu.
\end{equation}
$\mathfrak{D}^{n,k}_{\lambda,\mu}$ is a $\mathcal{K}(n)$-module for
the natural action:
\begin{equation*}
{X_F}\cdot A=\mathfrak{L}^{\mu}_{X_F}\circ A-(-1)^{|A||F|} A\circ
\mathfrak{L}^{\lambda}_{X_F}.
\end{equation*}
Thus, we have a filtration:
\begin{equation}\label{graded1}
\mathfrak{D}^{n,0}_{\lambda,\mu}\subset\mathfrak{D}^{n,\frac{1}{2}}_{\lambda,\mu}\subset
\mathfrak{D}^{n,1}_{\lambda,\mu}\subset
\mathfrak{D}^{n,\frac{3}{2}}_{\lambda,\mu}
\subset\cdots\subset\mathfrak{D}^{n,\ell-\frac{1}{2}}_{\lambda,\mu}\subset\mathfrak{D}^{n,\ell}_{\lambda,\mu}
\cdots
\end{equation}
Now, let us consider the density space
$\mathfrak{F}^n_\frac{2-n}{2}$ over the supercircle $S^{1|n}$. The
Berizin integral (\cite{NS, B, L1})
$\mathcal{B}_n:\mathfrak{F}^n_\frac{2-n}{2}\rightarrow\mathbb{R}$
can be given, for any $\varphi=\sum
f_{i_1,\dots,i_n}(x)\theta_1^{i_1}\cdots\theta_n^{i_n}\alpha_n^\frac{2-n}{2}$,
by the formula
$$\mathcal{B}_n(\varphi)=\int_{S^1} f_{1,\dots,1}dx.$$
\begin{proposition}
\label{conju} The Berizin integral $\mathcal{B}_n$ is
$\cK(n)$-invariant. That is, for any
$\varphi\in\mathfrak{F}^n_\frac{2-n}{2}$ and for any $H\in
C^\infty(S^{1|n})$, we have
$\mathcal{B}_n\left(\mathfrak{L}_{X_H}^\frac{2-n}{2}(\varphi)\right)=0$.
The product of densities composed with $\mathcal{B}_n$ yields a
bilinear $\cK (n)$-invariant form:
$$
\langle
~.~,.~\rangle:\mathfrak{F}^n_\lambda\otimes\mathfrak{F}^n_\mu\rightarrow\mathbb{R},\quad
\lambda+\mu=\frac{2-n}{2}.
$$
\end{proposition}
\begin{proofname}.
Note that $C^\infty(S)$ is assumed to be $\{f\in
C^\infty(\mathbb{R})~|~ f\text{ is } 2\pi\text{-periodic}\}$. For
$n=0$, we have $\mathfrak{L}_{X_H}^1(F)=HF'+H'F=(HF)'$, therefore,
$\mathcal{B}_0(\mathfrak{L}_{X_H}^1(F\alpha_n^1))=0$. For $n=1$,
using  equation (\ref{d-action}) for $\lambda=\frac{1}{2}$, we
easily show that
$\mathcal{B}_1(\mathfrak{L}_{X_H}^\frac{1}{2}(F\alpha_n^\frac{1}{2}))=0$.

 Let us consider $F=F_1+F_2\theta_n\in C^\infty(S^{1|n})$ and $H\in C^\infty(S^{1|n})$,
 where $\partial_n F_1=\partial_nF_2=\partial_nH=0$ with $\partial_i:=\frac{\partial}{\partial\theta_i}$. We easily prove that
$$
\mathfrak{L}^\frac{2-n}{2}_{X_H}F=\mathfrak{L}^\frac{2-n}{2}_{X_H}F_1+
\left(\mathfrak{L}^\frac{2-(n-1)}{2}_{X_H}F_2\right)\theta_n.
$$
So, we have
$\mathcal{B}_n\left(\mathfrak{L}_{X_H}^\frac{2-n}{2}\left(F\alpha_n^\frac{2-n}{2}\right)\right)=
\mathcal{B}_{n-1}\left(\mathfrak{L}_{X_H}^\frac{2-(n-1)}{2}\left(F_1\alpha_{n-1}^\frac{2-(n-1)}{2}
\right)\right)=0$.

On the other hand, by a direct computation, we show that
$$\mathfrak{L}_{X_{H\theta_n}}^\frac{2-n}{2}(F)=
(-1)^{|H|}\theta_n\left(\mathfrak{L}_{X_{H}}^\frac{2-(n-1)}{2}(F_1)
-\frac{1}{2}(H'F_1+HF'_1)\right)+\frac{1}{2}(-1)^{|F_2|}HF_2.$$ So,
it is clear that
$\mathcal{B}_n\left(\mathfrak{L}_{X_{H\theta_n}}^\frac{2-n}{2}\left(F\alpha_n^\frac{2-n}{2}\right)\right)=0$.
This completes the proof. \hfill$\Box$
\end{proofname}
\begin{cor}
 There exists a $\cK(n)$-invariant conjugation map:
$$*:\mathfrak{D}^{n,k}_{\lambda,\mu}\rightarrow
\mathfrak{D}^{n,k}_{\frac{2-n}{2}-\mu,\frac{2-n}{2}-\lambda}\quad\hbox{
defined by}\quad
\langle A\varphi,\psi\rangle=(-1)^{|A||\varphi|}\langle
\varphi,A^*\psi\rangle.
$$
Moreover, $*$ is $\cK(n)$-isomorphism
$\mathfrak{D}^{n,k}_{\lambda,\mu}\simeq\mathfrak{D}^{n,k}_{\frac{2-n}{2}-\mu,\frac{2-n}{2}-\lambda}$
for every $k\in{1\over 2}\mathbb{N}$.
\end{cor}
\section*{Acknowledgments}
We would like to  thank Valentin Ovsienko for many stimulating
discussions and valuable correspondences. Special thanks are due to
Charles H. Conley for his interest in this work and a number of
suggestions that have greatly improved this paper.


\end{document}